\documentclass[a4paper,10pt,reqno]{amsart}
\usepackage{amssymb}
\usepackage{color}
\usepackage[hyphens]{url}
\usepackage[square,numbers,comma,sort&compress]{natbib}
\usepackage{enumerate}
\usepackage{colortbl}
\usepackage{booktabs}
\usepackage{iba-algo}
\usepackage{graphics}
\usepackage{epsfig}
\usepackage{subfigure}
\usepackage{booktabs}
\usepackage{mathrsfs}   

\def\ve#1{\mathchoice{\mbox{\boldmath$\displaystyle\bf#1$}}
{\mbox{\boldmath$\textstyle\bf#1$}}
{\mbox{\boldmath$\scriptstyle\bf#1$}}
{\mbox{\boldmath$\scriptscriptstyle\bf#1$}}}

\newcommand\Z{\mathbb Z}   
\newcommand\R{\mathbb R}   
\newcommand\Po{\mathcal P}   
\newcommand\Cg{\mathcal C}   
\newcommand\In{\mathcal I}  
\newcommand\B{\mathcal B}   
\newcommand\E{\mathcal E}   

\DeclareMathOperator{\init}{in}

\DeclareMathOperator{\conv}{conv}   
\DeclareMathOperator{\rank}{rank}

\let\epsilon=\varepsilon
\usepackage{ifthen}
\makeatletter
\newcommand{\DeclareBracket}[3]{
  \newcommand{#1}[2][]{%
  \ifthenelse%
  {\equal{##1}{}}%
  {\left#2##2\right#3}%
  {\csname ##1l\endcsname#2##2\csname ##1r\endcsname#3}}}    
\makeatother
\DeclareBracket\abs||           
\DeclareBracket\norm\|\|        
\DeclareBracket\floor\lfloor\rfloor
\DeclareBracket\ceil\lceil\rceil
\DeclareBracket\set\{\}
\DeclareBracket\paren()
\DeclareBracket\inner\langle\rangle
\DeclareBracket\fractional\{\}
 
\usepackage{amsthm}
\newtheorem{theorem}{Theorem}%
\makeatletter
\newtheorem{lemma}{Lemma}
\renewcommand*{\c@lemma}{\c@theorem}
\renewcommand*{\p@lemma}{\p@theorem}

\renewcommand*{\c@conjecture}{\c@theorem}
\renewcommand*{\p@conjecture}{\p@theorem}

\newtheorem{proposition}{Proposition}
\renewcommand*{\c@proposition}{\c@theorem}
\renewcommand*{\p@proposition}{\p@theorem}

\newtheorem{corollary}{Corollary}
\renewcommand*{\c@corollary}{\c@theorem}
\renewcommand*{\p@corollary}{\p@theorem}

\renewcommand*{\c@observation}{\c@theorem}
\renewcommand*{\p@observation}{\p@theorem}

\theoremstyle{definition}

\renewcommand*{\c@problem}{\c@theorem}
\renewcommand*{\p@problem}{\p@theorem}

\newtheorem{definition}{Definition}
\renewcommand*{\c@definition}{\c@theorem}
\renewcommand*{\p@definition}{\p@theorem}

\newtheorem{remark}{Remark}
\renewcommand*{\c@remark}{\c@theorem}
\renewcommand*{\p@remark}{\p@theorem}

\newtheorem{example}{Example}
\renewcommand*{\c@example}{\c@theorem}
\renewcommand*{\p@example}{\p@theorem}

\newtheorem{algorithm}{Algorithm}
\renewcommand*{\c@algorithm}{\c@theorem}
\renewcommand*{\p@algorithm}{\p@theorem}

\makeatother
\usepackage{hyperref}

\title{Volumes and Tangent Cones of Matroid Polytopes}

\author{David C. Haws}
\address{David C. Haws: Department of Statistics,
University of Kentucky,
Lexington, KY 40508,
Copyright 2011}
\email{\texttt dchaws@gmail.com}

\begin{document}

\begin{abstract}
De Loera et al.\ 2009, showed that when the rank is fixed the Ehrhart
polynomial of a matroid polytope can be computed in polynomial time when the
number of elements varies. A key to proving this is the fact that the number
of simplicial cones in any triangulation of a tangent cone is bounded
polynomially in the number of elements when the rank is fixed. The authors
speculated whether or not the Ehrhart polynomial could be computed in
polynomial time in terms of the number of bases, where the number of elements
and rank are allowed to vary. We show here that for the uniform matroid of rank
$r$ on $n$ elements, the number of simplicial cones in any triangulation of a
tangent cone is $n-2 \choose r-1$. Therefore, if the rank is allowed to vary,
the number of simplicial cones grows exponentially in $n$. Thus, it is unlikely
that a Brion-Lawrence type of approach, such as Barvinok's Algorithm, can compute the
Ehrhart polynomial efficiently when the rank varies with the number of
elements. To prove this result, we provide a triangulation in which the maximal
simplicies are in bijection with the spanning thrackles of the complete
bipartite graph $K_{r,n-r}$.
\end{abstract}

\maketitle

\section{Introduction}
  
Recall that a \emph{matroid} $M$ is a finite
collection $\mathcal{F}$ of subsets of $[n] = \{1,2,\dots,n\}$ called
\emph{independent sets}, such that the following properties are satisfied:
{\bf (1)} $\emptyset \in \mathcal{F}$, {\bf (2)} if $X \in \mathcal{F}$ and $Y \subseteq X$ then $Y \in \mathcal{F}$, {\bf (3)} if $U, V \in \mathcal{F}$ and $|U| = |V| + 1$ there exists $x \in U \setminus V$ such that $V \cup x \in \mathcal{F}$. In this paper
we investigate  convex polyhedra associated with matroids.

%

Similarly, recall that a matroid $M$ can be
defined by its \emph{bases}, which are the inclusion-maximal independent sets.
The bases of a matroid $M$ can be recovered by its rank function $\varphi$. For
the reader we recommend \cite{Oxley1992Matroid-Theory} or
\cite{Welsh1976Matroid-Theory} for excellent introductions to the theory of
matroids.


Now we introduce the main object of this paper. Let $\mathcal{B}$ be the set
of bases of a matroid $M$. If $B = \{\sigma_1,\ldots,\sigma_r\} \in
\mathcal{B}$, we define the \emph{incidence vector of B} as $\ve  e_B :=
\sum_{i=1}^r \ve e_{\sigma_i}$, where $\ve e_j$ is the standard elementary
$j$th vector in $\R^n$. The \emph{matroid polytope} of $M$ is defined as
$\Po(M) := \conv \{\, \ve e_B  \mid B \in \mathcal{B} \, \}$, where
$\conv(\cdot)$ denotes the convex hull. This is different from the well-known
\emph{independence matroid polytope}, $\Po^{\In} (M) := \conv \{ \, \ve e_I
\mid  I \subseteq B \in \mathcal B \, \}$, the convex hull of the incidence
vectors of all the independent sets. We can see that $\Po(M) \subseteq
\Po^\In(M)$ and $\Po(M)$ is a face of $\Po^\In(M)$ lying in the hyperplane
$\sum_{i=1}^n x_i = \rank(M)$, where $\rank(M)$ is the cardinality of any basis
of $M$.

Recall that given an integer $k > 0$ and a polytope $\Po \subseteq \R^n$ we
define $k \Po := \{\, k \ve \alpha  \mid  \ve \alpha \in \Po \, \}$ and the
function $i(\Po,k) := \#(k\Po \cap \Z^n) $, where we define $i(\Po,0) :=1$. It
is well known that for integral polytopes, as in the case of matroid polytopes,
$i(\Po,k)$ is a polynomial, called the \emph{Ehrhart polynomial} of $\Po$.
Moreover the leading coefficient of the Ehrhart polynomial is the
\emph{normalized volume} of $\Po$, where a unit is the volume of the
fundamental domain of the affine lattice spanned by $\Po$
\cite{Stanley1996Combinatorics-a}.  In \cite{De-Loera:2009uq} the following was
shown:

\begin{theorem}[ Theorem 1 in \cite{De-Loera:2009uq}]
\label{thm:volume}
  Let $r$ be a fixed integer.
  Then there exist algorithms whose input data consists of a number $n$ and an evaluation oracle for 
  \begin{enumerate}[\rm(a)]
  \item a rank function~$\varphi$ of a matroid~$M$ on $n$ elements
    satisfying $\varphi(A) \leq r$ for all~$A$, or
  \item an integral polymatroid rank function~$\psi$ satisfying $\psi(A) \leq r$ for
    all~$A$,
  \end{enumerate}
  that compute in time polynomial in~$n$ the Ehrhart polynomial (in
  particular, the volume) of the matroid polytope $\Po(M)$, the independence
  matroid polytope $\Po^\In(M)$, and the polymatroid~$\Po(\psi)$,
  respectively.
\end{theorem}

The proof of \autoref{thm:volume} relied on four important facts when the
rank is fixed: (1) The number of bases is polynomially bounded, (2) every
triangulation of a tangent cone of the matroid polytope has a polynomial number
of maximal simplicial cones, (3) a triangulation of a tangent cone can be done
in polynomial time, and (4) every triangulation of a tangent cone of the
matroid polytope is unimodular. The first item follows easily from the rank
being fixed, implying there are at most $n \choose r$ bases. The third item is
relatively straightforward using the pulling triangulation, given item two
holds. The proof of item two in \cite{De-Loera:2009uq} (Lemma 10) used a bound
on the volume of the subpolytope given by a vertex and all its adjacent
vertices. Item four does not rely on the rank being fixed at all.

The authors of \cite{De-Loera:2009uq} speculated whether or not the Ehrhart
polynomial of a matroid polytope could be computed in polynomial time with
respect to the number of basis, regardless of the rank. The primary limitation
to proving this result seemed to be item two above. However, we show the
following: 

\begin{theorem}
\label{thm:main}
Let $U^{r,n}$ be the uniform matroid of rank $r$ with $n$ elements. There
are $n-2 \choose r-1$ simplicial cones in any triangulation of a tangent cone of
the matroid base polytope of $U^{r,n}$.
\end{theorem}

Thus, it is unlikely that the Erhrart polynomial of a matroid base polytope can
be computed in polynomial time when the number of ground elements $n$ varies and the
rank is not fixed.  That is, a Brion-Lawrence type of approach (Barvinok's Algorithm \cite{bar})
to compute the Ehrhart polynomial is likely not computationally
efficient. If the rank is allowed to vary, \autoref{thm:main} states that the
number of simplicial cones in the triangulation of any tangent cone of
$U^{n/2,n}$, where $n$ even, is $n-2 \choose n/2 -1$, the central
binomial coefficient. And, it is known that $n-2 \choose n/2 -1 $ grows exponentially
in $n$.

\section{Gr\"obner Bases and Triangulations}
Notation and ideas for many of the proofs in this section are taken from
\cite{Sturmfels1996Grobner-Bases-a} (which is in turn was drawn from
\cite{Loera:1995fk}), which covers Gr\"obner bases and triangulations.

The edges of the matroid polytope have the following important property.
\begin{lemma}[See Theorem 4.1 in \cite{Gelfand1987Combinatorial-g}, Theorem 5.1 and Corollary 5.5 in \cite{Topkis1984Adjacency-on-Po}] \label{lem:adj}
    Let $M$ be a matroid. 
    \begin{itemize}
    \item[A)]   Two vertices $\ve e_{B_1}$ and $\ve e_{B_2}$ are adjacent in $\Po(M)$ if and only if $\ve e_{B_1} - \ve e_{B_2} = \ve e_i - \ve e_j$ for some $i,j$.
    \item[B)]   If two vertices $\ve e_{I_1}$ and $\ve e_{I_2}$ are adjacent
      in $\Po^\In (M)$ then $\ve e_{I_1} - \ve e_{I_2} \in \{ \, \ve e_i - \ve
      e_j, \, \ve e_i, \, -\ve e_j \, \}$ for some $i,j$. Moreover if $\ve v$
      is a vertex of $\Po^\In (M)$ then all adjacent vertices of $\ve v$ can
      be computed in polynomial time in $n$, even if the matroid $M$ is only
      presented by an evaluation oracle of its rank function~$\varphi$.
    \end{itemize}
\end{lemma}

In this section we study the matroid polytopes of uniform matroids. Let
$U^{r,n}$ denote the uniform matroid of rank $r$ on $n$ elements. If $B$ is a
basis of $U^{r,n}$, then $B\backslash \{i\} \cup \{j\}$ is an adjacent basis on
$\Po(U^{r,n})$, for all $i \in B$ and $j \in [n] \backslash B$. Since we study
the uniform matroid, without loss of generality, we focus on the basis
$B = \{1,\ldots,r\}$ with incidence vector 
\begin{equation*}
[\underbrace{1,\ldots,1}_{r},\underbrace{0,\ldots,0}_{n-r}].
\end{equation*}
The bases adjacent to $B$ are then $\{1,\ldots,r\} \backslash \{i\} \cup \{j\}$
where $1\leq i \leq r$ and $r+1 \leq j \leq n$. Shifting by the vector $\ve
e_B$, we are interested in the pointed cone $C$ with rays $\ve e_i - \ve e_j$,
for all $1\leq i \leq r$ and $r+1 \leq j \leq n$. We define
\begin{equation*}
\E_{r,n} := \{\, \ve e_j - \ve e_i \mid 1\leq i \leq r,\, r+1 \leq j \leq n \,\}.
\end{equation*}

Note that the points $\E_{r,n}$ lie in the hyperplane
$[\underbrace{1,\ldots,1}_{r},\underbrace{-1,\ldots,-1}_{n-r}] \cdot \ve x =
0$, which does not contain the origin. We are interested in triangulations of
$C$ into simplicial cones. Note that any triangulation of $C$ is in
bijection with triangulations of the points $\E_{r,n}$. Hence, we study
triangulations of $\E_{r,n}$. To simplify matters and better relate to material
in \cite{Sturmfels1996Grobner-Bases-a} we will focus on triangulations of 

\begin{equation*}
\B_{r,n} := \{\, \ve e_i + \ve e_j \mid 1\leq i \leq r,\, r+1 \leq j \leq n \,\}.
\end{equation*}

Note that $\conv(\E_{r,n})$ can be mapped to $\conv(\B_{r,n})$ by the unimodular involution

\begin{equation*}
\left( \begin{array}{cc} -I_r & 0 \\ 
                        0 & I_{n-r} \end{array} \right).
\end{equation*}
Therefore triangulations of $\conv(\E_{r,n})$ and $\conv(\B_{r,n})$ are in
bijection and their volumes are equal. Note that $\B_{r,n}$ is a sub-polytope
of the \emph{second hypersimplex} 
\begin{equation*}
\mathcal A_n := \{\, \ve e_i + \ve e_j \mid 1 \leq i < j \leq n \, \}.
\end{equation*}

We now follow closely the notation and proofs of Chapter 9 of
\cite{Sturmfels1996Grobner-Bases-a}.  In it, Sturmfels gives a
unimodular triangulation of $\mathcal A_n$ into $2^{n-1}-n$ simplices.

\begin{proposition}
The dimension of $\conv(\B_{r,n})$ is $n-2$.
\end{proposition}
\begin{proof}
The points $\B_{r,n}$ lie in the hyperplanes $[1,\ldots,1]\cdot \ve x = 2$ and 
$[\underbrace{1,\ldots,1}_{r},\underbrace{-1,\ldots,-1}_{n-r}] \cdot \ve x =
0$. It is not difficult to see there are $n-2$ linearly independent vectors in 
$\B_{r,n}$.
\end{proof}

\begin{remark}
The set of column vectors $\B_{r,n}$ is the vertex-edge incidence matrix of
the complete bipartite graph $K_{r,n-r}$.
\end{remark}

The toric ideal $I_{\B_{r,n}}$ is the kernel of the map

\begin{equation*}
\Phi : k[ x_{ij} \, : \, 1 \leq i \leq r, \, r+1 \leq j \leq n] \rightarrow k[t_1,\ldots,t_n], \, x_{ij} \mapsto t_it_j.
\end{equation*}

The variables $x_{ij}$ are indexed by the edges of the complete bipartite graph
$K_{r,n-r}$. We identify the vertices of $K_{r,n-r}$ with the vertices of a
planar embedding of the complete bipartite graph on $r$ and $n-r$ vertices.  By
an \emph{edge} we mean a closed line segment between two vertices in the
complete bipartite graph on $r$ and $n-r$ vertices. The \emph{weight} of the
variable $x_{ij}$ is the number of edges of $K_{r,n-r}$ which do not meet the
edge $(i,j)$.

\begin{figure}
\includegraphics[height=1in]{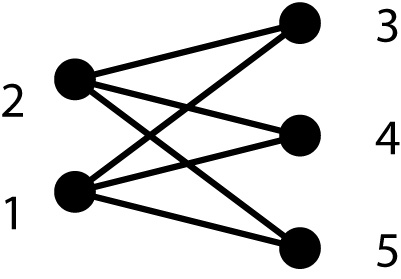}
\caption{Complete bipartite graph $K_{2,3}$.}
\label{fig:K2_3}
\end{figure}

\begin{example}
For $K_{2,3}$ in \autoref{fig:K2_3}, variables $x_{13},x_{25}$ have weight $0$,
variable $x_{14},x_{24}$ have weight $1$, and variables $x_{15}, x_{23}$ have weight $2$.
\end{example}

\begin{remark}
Throughout this section, we will draw the complete bipartite graph $K_{s,t}$ as
in \autoref{fig:K2_3}: The $s$ vertices $\{1,\ldots,s\}$ are drawn vertically
on the left and labeled from bottom to top and the $t$ vertices
$\{s+1,\ldots,s+t\}$ are drawn vertically on the right and labeled from top to
bottom. This is done to match closely with \cite{Sturmfels1996Grobner-Bases-a},
where the complete graphs $K_n$ are drawn on the $n$-gon with labels in
clockwise order. Thus, for the drawings of the complete bipartite graph
$K_{s,t}$ we can also talk about the $(s+t)$-gon given by the $s+t$ points.
\label{rem:draw}
\end{remark}

Given any pair of non-intersecting edges $(i,j),(k,l)$ of $K_{r,n-r}$ the pair
$(i,l),(j,k)$ meets in a point (intersect). With disjoint edges $(i,j),(k,l)$
we associate the binomial $x_{ij}x_{kl} - x_{il}x_{jk}$. We denote by $\Cg$ the
set of all binomials obtained in this fashion, and by $in_\succ (\Cg)$ the set
of their initial monomials.  Here $\succ$ denotes the term order that refines
the partial order on monomials specified by these weights.

\begin{example}
Let $r=2$, $n=5$. Then $\Cg$ is
\begin{equation*}
\big\{ \, \underline{x_{13}x_{25}}-x_{15}x_{23},\underline{x_{13}x_{24}}-x_{14}x_{23},\underline{x_{14}x_{25}}-x_{15}x_{24} \big\}.
\end{equation*}
\end{example}

\begin{theorem}
The set $\Cg$ is the reduced Gr\"obner basis of $I_{B_{r,n}}$ with respect to $\succ$.
\label{thm:grobner}
\end{theorem}

The proof of \autoref{thm:grobner} follows nearly verbatim of the proof of
Theorem 9.1 in \cite{Sturmfels1996Grobner-Bases-a}, with only minor
modification to handle the complete bipartite graph. The proof of
\autoref{thm:grobner} is given in the Appendix.

\begin{remark}[Remark 9.2 in \cite{Sturmfels1996Grobner-Bases-a}]
The set $\Cg$ is the reduced Gr\"obner basis for $I_{B_{r,n}}$ with respect to
the purely lexicigraphic term order induced by the following variable ordering:
\begin{equation*}
x_{ij} \succ x_{kl} \quad \text{if and only if} \quad i < k \text{ or } (i=k \text{ and } j > l).
\end{equation*}
\begin{proof}
For any ordered quadruple $1 \leq i < j < k < l \leq n$, the intersecting pair
of edges is $\{(i,k),(j,l)\}$. We must show that the monomial $x_{ik}x_{jl}$ is
smaller than both $x_{ij}x_{jk}$ in the given term order. But this holds since
$x_{kl} \succ x_{jk} \succ x_{jl} \succ x_{ij} \succ x_{ik} \succ x_{il}$.
\end{proof}
\end{remark}

Following identical logic in Chaption 9 of \cite{Sturmfels1996Grobner-Bases-a},
we apply \autoref{thm:grobner} to give an explicit triangulation and
determine the normalized volume of $\conv(B_{r,n})$. By Theorem 8.3 in
\cite{Sturmfels1996Grobner-Bases-a}, the square-free monomial ideal $\langle
\init_\succ (\Cg) \rangle = in_\succ (I_{B_{r,n}})$ is the Stanley-Reisner
ideal of a regular triangulation $\Delta_\succ$ of $\conv(B_{r,n})$. The
simplices in $\Delta_\succ$ are the supports of the standard monomials. All
maximal simplicies in $\Delta_\succ$ have unit normalized volume by Corollary
8.9 in \cite{Sturmfels1996Grobner-Bases-a}, Corollary 63 in \cite{Haws:2009uq},
or Lemma 8 in \cite{De-Loera:2009uq}. We observed before that elements of
$\init_\succ(\Cg)$, the minimally non-standard monomials, are supported
on pairs of disjoint edges.

\begin{corollary}
\label{cor:maxsimp}
The simplices of $\Delta_\succ$ are the subgraphs of $K_{r,n-r}$ with the
property that any pair of edges intersects in the convex embedding of the graph
given in \autoref{rem:draw}.
\end{corollary}

Now we identify subgraphs of $K_{r,n-r}$ with subpolytopes of $\conv(B_{r,n})$:
A subgraph $H$ is identified with the convex hull of the column vectors of its
vertex-edge incidence matrix.

\begin{definition}
Let $G$ be a graph with edge set $E$ embedded in the plane. Recall the edges
$\widehat E \subseteq E$ are a \emph{thrackle} if every pair of edges
intersects.
\end{definition}

\begin{proposition}
\label{prop:simp}
A subpolytope $\sigma$ of $\conv(B_{r,n})$ is a $(n-2)$-dimensional simplex if
and only if the coresponding subgraph $H$ is a spanning tree of $K_{r,n-r}$.
The normalized volume of $\sigma$ is $1$.
\end{proposition}
\begin{proof}
Suppose $H$ supports a $(n-2)$-simplex. Let $M_H$ be the $\{0,1\}$-incidence
matrix of $H$. This matrix is non-singular which implies it is spanning and all
cycles are odd (if any). But, the complete bipartite graph does not have odd
cycles. Therefore $H$ is acyclic and spanning. Any acyclic spanning subgraph of
$K_{r,n-r}$ with $n-1$ edges is a spanning tree.

Conversely, if $H$ is a spanning tree of $K_{r,n-r}$, it contains $n-1$ edges,
is acyclic and its incidence matrix is non-singular implying it is a
$(d-2)$-simplex. 

There exists some vertex $v \in [n]$ such that the degree of $v$ is $1$.
Performing cofactor expansion on the $v$th row we see that $M_H = M_{H-v}$.
Repeating we see that $|M_H| = 1$.
\end{proof}

\begin{theorem}
\label{thm:maxsimpspanthrackles}
The maximal simplices of the triangulation $\Delta_\succ$ are the spanning trees
of $K_{r,n-r}$ with $n-1$ edges and with the property that any pair of edges
intersects.
\end{theorem}
\begin{proof}
Follows from \autoref{cor:maxsimp} and \autoref{prop:simp}.
\end{proof}

\begin{proposition}
Every thrackle of $K_{s,t}$ is acyclic.
\end{proposition}
\begin{proof}
All cycles in $K_{s,t}$ must be even, but as shown in the proof of 
\autoref{thm:grobner}, any even cycle will contain two edges that do not cross.
\end{proof}

\autoref{thm:maxsimpspanthrackles} states that the simplices of
$\Delta_\succ$ are the spanning thrackles of $K_{r,n-r}$. It is known that the
number of spanning trees of $K_{r,n-r}$ is $r^{n-r-1}n^{r-1}$, but from Lemma
10 in \cite{De-Loera:2009uq} and by observation, not all spanning trees are
thrackles. Below we prove that the number of spanning thrackles is simply
a binomial coefficient.

We first observe a basic fact about thrackles on $K_{r,n-r}$.
\begin{proposition}
\label{prop:thracklecontseg}
Let $K_{s,t}$ be the complete bipartite graph with the planar embedding in
\autoref{rem:draw}, $H$ a thrackle of $K_{s,t}$, $1 \leq i \leq s$, $s+1 \leq j_1 < j_2
\leq s+t$, $j_1 +1 < j_2$ and $\{(i,j_1),(i,j_2)\}$ edges of $H$. Then for all
$k \in [s] \backslash \{i\}$ and $j_1 < l < j_2$, $(k,l)$ does not intersect
every edge of $H$.
\end{proposition}
\begin{proof}
Let $K_{s,t}$ be the complete bipartite graph with the planar embedding
in \autoref{rem:draw}, $H$ a thrackle of $K_{s,t}$, $1 \leq i \leq s$, $s+1 \leq j_1 \leq j_2
\leq s+t$, $j_1 +1 < j_2$ and $\{(i,j_1),(i,j_2)\}$ edges of $H$.  If $k > i$
then $(k,l)$ does not intersect $(i,j_2)$. If $k < i$ then $(k,l)$ does not
intersect $(i,j_1)$.
\end{proof}

\begin{figure}
\includegraphics[height=1in]{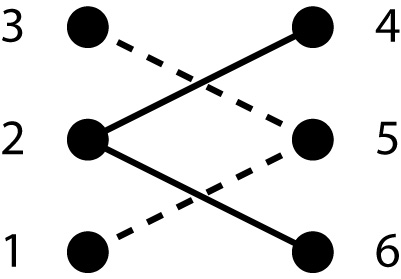}
\caption{A thrackle (solid) with edges $\{(2,4),(2,6)\}$. No other edges
(dashed) incident to vertex $5$ and not incident to vertex $2$ can be added to
the thrackle to give a new thrackle.}
\label{fig:thracklecontsegex}
\end{figure}

For an example of \autoref{prop:thracklecontseg}, see 
\autoref{fig:thracklecontsegex}. From the above proposition we get the following
corollary
\begin{corollary}
\label{cor:spanthracklecontseg}
Let $K_{s,t}$ be the complete bipartite graph with the planar embedding in
\autoref{rem:draw}. If $H$ is a spanning thrackle of $K_{s,t}$, $1 \leq i
\leq s$, $s+1 \leq j_1 \leq j_2 \leq s+t$, $j_1 +1 < j_2$ and
$\{(i,j_1),(i,j_2)\}$ edges of $H$, then for all $j_1 \leq j_3 \leq j_2$,
$(i,j_3)$ is an edge of $H$.
\end{corollary}

\begin{proposition}
\label{prop:firstedge}
Let $K_{s,t}$ be the complete bipartite graph with the planar embedding in
\autoref{rem:draw} and $H$ a spanning thrackle of $K_{s,t}$. Then
$\{1,s+1\}$ must be in $H$.
\end{proposition}
\begin{proof}
Since $H$ is spanning, some edge of $H$ must be incident to $s+1$. But, this
edge will not intersect any edge $\{1,j\}$ where $s+1 < j \leq s+t$.
\end{proof}

\begin{definition}
Let $K_{s,t}$ be the complete bipartite graph with the planar embedding in
\autoref{rem:draw}. We define $f(s,t)$ to be number of spanning thrackles of
$K_{s,t}$.
\end{definition}

By observation we see that $f(1,1) = 1$ and $f(1,t)=1$ for all $1 \leq t$. Also
note that $f(s,t) = f(t,s)$.
\begin{proposition}
\label{prop:rec}
Let $K_{s,t}$ be the complete bipartite graph with the planar embedding in
\autoref{rem:draw}. The number of spanning thrackles such that the only
edges incident to vertex $1$ are $\{(1,s+1),\ldots,(1,i)\}$ where $s+1 \leq
i \leq t+s$ is $f(s-1,t+s+1-i))$.
\end{proposition}
\begin{proof}
Let $K_{s,t}$ be the complete bipartite graph with the planar embedding in
\autoref{rem:draw} and $H$ a spanning thrackle such that the only edges
incident to vertex $1$ are $\{(1,s+1),\ldots,(1,i)\}$ where $s+1 \leq i \leq
t+s$. For any $1 < j \leq s$ and $s+1 \leq k < i$ we have that $(j,k)$ does not
intersect the edges $\{(1,s+1),\ldots,(1,i)\}$. Therefore, the number of
spanning thrackles such that the only edges incident to vertex $1$ are
$\{(1,s+1),\ldots,(1,i)\}$ where $s+1 \leq i \leq t+s$ equals
$f(s-1,t+s+1-i))$.
\end{proof}

\begin{figure}
\includegraphics[height=0.9in]{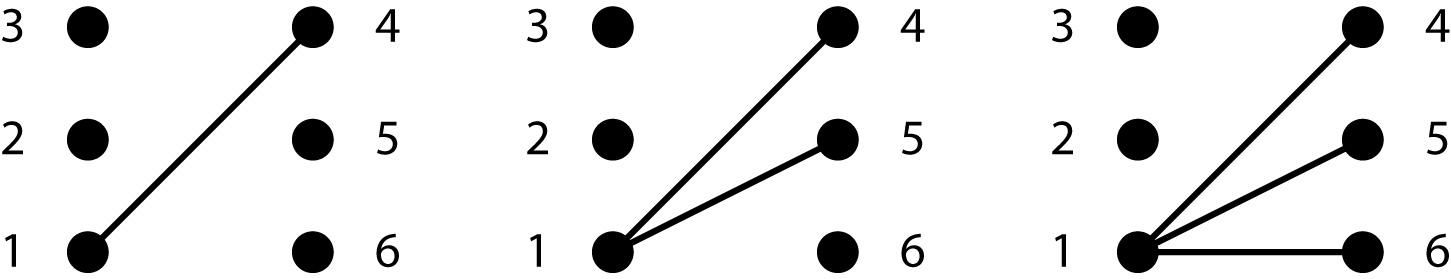}
\caption{For $K_{3,3}$, the spanning thrackles fall into three disjoint cases.
The only edge incident to $1$ is $(1,4)$ (left).  The only edge incident to $1$
is $\{(1,4),(1,5)\}$ (middle). The only edge incident to $1$ is
$\{(1,4),(1,5),(1,6)\}$ (right).}
\label{fig:thracklediscases}
\end{figure}

Using \autoref{prop:rec}, we find a recurrence relation
for the number of spanning thrackles of $K_{s,t}$ by dividing the spanning
thrackles $H$ into disjoint cases (see \autoref{fig:thracklediscases}):
\begin{enumerate}
\item The only edge incident to $1$ is $\{(1,s+1)\}$.
\item The only edges incident to $1$ are $\{(1,s+1),(1,s+2)\}$.
\item The only edges incident to $1$ are $\{(1,s+1),(1,s+2),(1,s+3)\}$.\\
$\vdots$
\item The only edges incident to $1$ are $\{(1,s+1),(1,s+2),\ldots,(1,s+t)\}$.
\end{enumerate}

In item (1), if the only edge incident to $1$ in $H$ is $\{(1,s+1)\}$, then the
number of spanning thrackles satisfying this condition is equal to
$f(s-1,s+t)$. Similarly, if the only edges incident to $1$ in $H$ are
$\{(1,s+1),(1,s+2)\}$, then the number of spanning thrackles satisfying this
condition is equal to $f(s-1,s+t-1)$.  This leads to the recursion
relation 
\begin{equation}
f(s,t) = \sum_{i=0}^{t-1} f(s-1,t-i).
\label{eq:thracklerecrel}
\end{equation}

\begin{lemma}
\label{lem:uniqpartthrackle}
Let $K_{s,t}$ be the complete bipartite graph with the planar embedding in
\autoref{rem:draw} and $H$ a spanning thrackle of $K_{s,t}$. The edges
\begin{equation}
\label{eq:uniqedges}
\left\{\, (i,j) \in H \mid  \nexists k,\, k<j,\, (i,k) \in H \, \right\}
\end{equation}
uniquely determine $H$.
\end{lemma}
\begin{figure}
\includegraphics[height=1.2in]{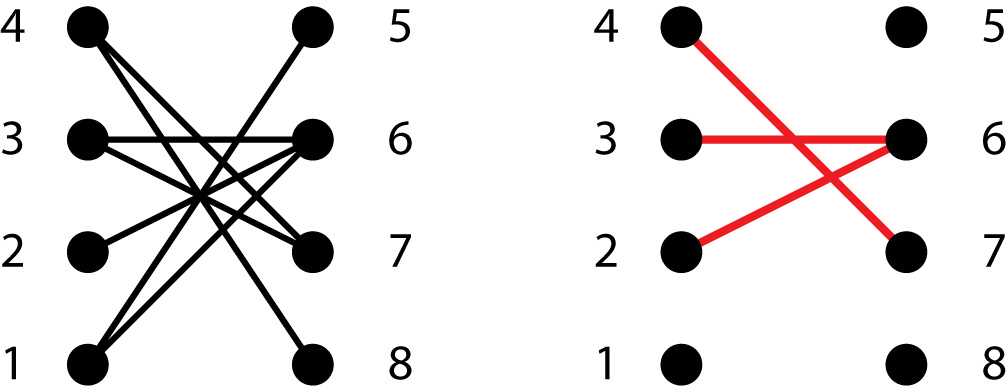}
\caption{For the spanning thrackle (left) of $K_{4,4}$ the red edges (right) uniquely determine the thrackle.}
\label{fig:thrackleuniqparts}
\end{figure}

See \autoref{fig:thrackleuniqparts} for an example.
\begin{proof}[Proof of \autoref{lem:uniqpartthrackle}]
By \autoref{cor:spanthracklecontseg} and \autoref{prop:firstedge}
the vertices incident to vertex $1$ must be an interval $[s+1,i_1]$ where $s+1
\leq i_1 \leq s+t$. Similarly, by the recursive argument in 
\autoref{prop:rec}, the vertices incident to vertex $2$ must be an interval
$[i_1,i_2]$ where $ i_1 \leq i_2 \leq s+t$. Continuing this argument, the
vertices $1,2,\ldots,s$ are incident to the interval of vertices
$[s+1,i_1],[i_1,i_2],[i_2,i_3],\ldots,[i_{s-1},s+t]$ respectively.
Thus, $H$ is uniquely determined by $i_1,\ldots,i_{s-1}$ and these
are precisely what are given in Equation \eqref{eq:uniqedges}.
\end{proof}

\begin{remark}
Note that if $(i,j)$ is an edge in the set \ref{eq:uniqedges}, then for all
$(k,l)$ in the set \ref{eq:uniqedges} such that $k > i$ we have $l \geq j$.
\end{remark}

Recall from \cite{Stanley1997Enumerative-Com} that a \emph{weak composition} of
a positive integer $p$ is an ordered sum of $q$ non-negative integers which
sums to $p$.  It is known that the number of such weak composition is $p+q-1
\choose q-1$.

\begin{corollary}
\label{cor:numthrackles}
The number of spanning thrackles $f(s,t)$ of the complete bipartite graph
$K_{s,t}$ embedded in the plane as in \autoref{rem:draw} is equal to $s+t-2
\choose s-1$.  I.e., 
\begin{equation}
\label{eq:numthrackles}
f(s,t) = {s+t-2 \choose s-1} = {s+t-2 \choose t-1}.
\end{equation}
\end{corollary}
\begin{proof}
\autoref{lem:uniqpartthrackle} states that the edges of Equation \eqref{eq:uniqedges}
uniquely determine the spanning thrackle. But, it can be seen that these edges
are in bijection with the weak compositions of $t-1$ into $s$ parts.  Hence
\autoref{eq:numthrackles} holdes.
\end{proof}

Now we are equipped to prove \autoref{thm:main}.

\begin{proof}[Proof of \autoref{thm:main}]
By \autoref{thm:maxsimpspanthrackles} the simplices of $\Delta_{\succ}$ are
the maximal spanning thrackles, and from \autoref{cor:numthrackles}, we
know this to be $n-2 \choose r-1$. Now $\Delta_\succ$ is but one possible
triangulation of the tangent cone. But, by Corollary 8.9 in
\cite{Sturmfels1996Grobner-Bases-a}, Corollary 63 in \cite{Haws:2009uq}, or
Lemma 8 in \cite{De-Loera:2009uq}, we know that any triangulation of a tangent
cone of a matroid polytope will be composed of unimodular cones. Hence it
will have the same normalized volume, namely $n-2 \choose r-1$.
\end{proof}

Finally, we offer an alternative method to prove Equation \eqref{eq:numthrackles}.  
We give a function $\Phi$ from the space spanning thrackles of $K_{s,t}$
embedded as in \autoref{rem:draw} to the space of $\{0,1\}$ sequences with
exactly $s-1$ zeros and $t-1$ ones. Let $H$ be a spanning thrackle. The string
$\Phi(H)$ is given by the algorithm:

\begin{algorithm}[]
    \mbox{}
    \vskip .25cm
    \label{alg:thracklebij}
    \begin{center}
        \begin{tabular}{l}
            \begin{minipage}{.85\linewidth}
                \begin{algorithmic}
                    \item[Input:] A spanning thrackle $H$ of $K_{s,t}$ embedded as in \autoref{rem:draw}.
                    \item[Output:] A $\{0,1\}$ string with exactly $s-1$ zeros and $t-1$ ones.
                    \STATE Let $S$ = ``'';
                    \FOR{every edge $v \in [s+1,s+t-1]$} 
                        \FOR{every vertex $w \neq 1$ adjacent to $v$ in $H$ and $w$ is unmarked} 
                            \STATE S = S + ``0''.
                            \STATE Mark $w$.
                        \ENDFOR
                        \STATE S = S + ``1''.
                    \ENDFOR
                    \FOR{every vertex $w \neq 1$ adjacent to $s+t$ in $H$ and $w$ is unmarked} 
                        \STATE S = S + ``0''.
                        \STATE Mark $w$.
                    \ENDFOR
                \end{algorithmic}
            \end{minipage}\\
        \end{tabular}
    \end{center}
\end{algorithm}

\begin{figure}
\includegraphics[height=2.5in]{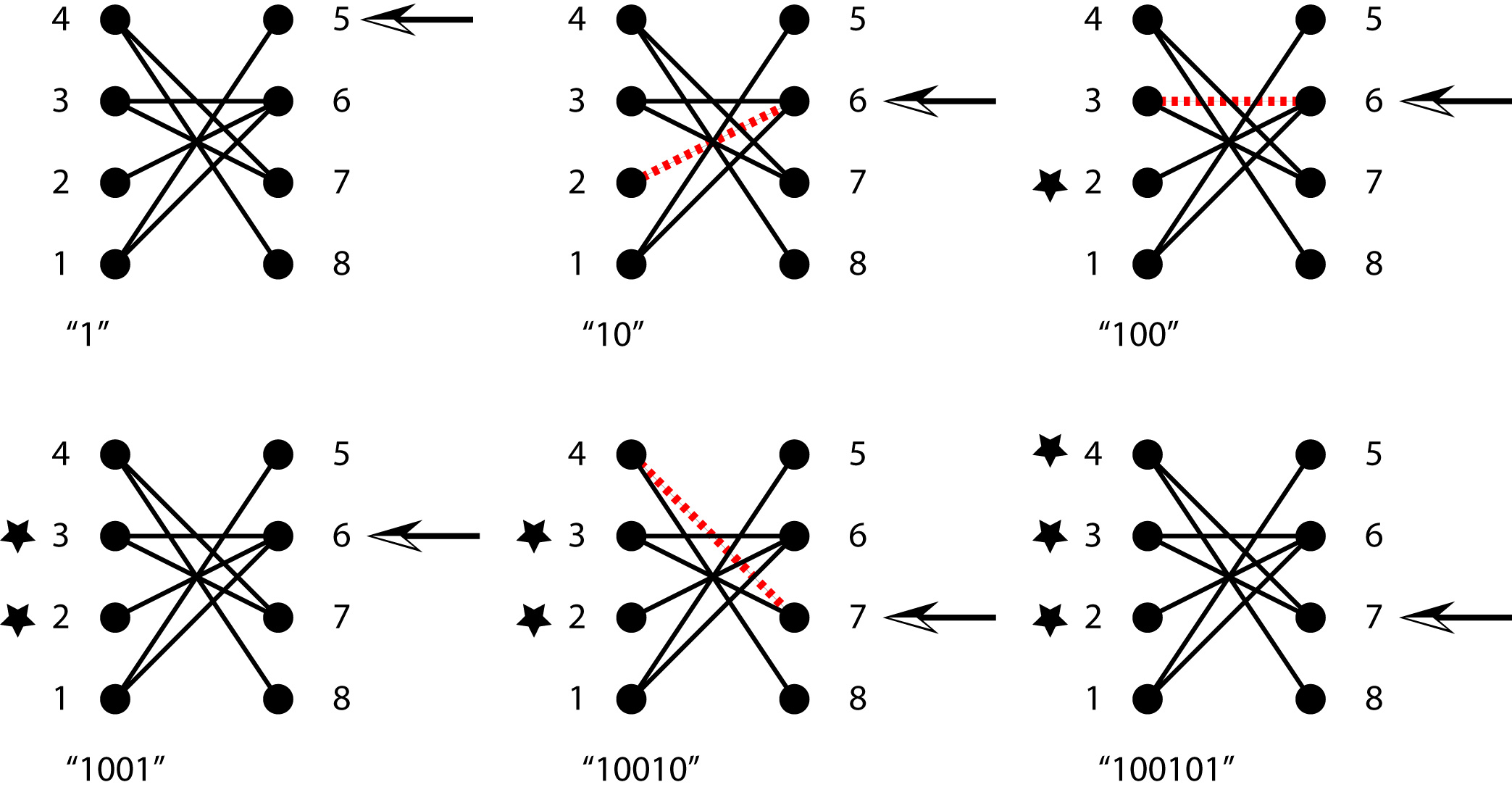}
\caption{An example of \autoref{alg:thracklebij} mapping of a spanning
thrackle to a $\{0,1\}$ sequence with $s-1$ zeros and $t-1$ ones.}
\label{fig:thracklebij}
\end{figure}

Following this, it is known the number of such $\{0,1\}$ strings is $s+t-2
\choose s-1$.  For an example of \autoref{alg:thracklebij} (i.e. $\Phi$), see
\autoref{fig:thracklebij}. 


\section{Discussion}
If $M$ is not the uniform matroid, then the extreme rays of a tangent cone will
not necessarily be $\mathcal E_{r,n}$. However, the arguments in 
\autoref{thm:grobner} will still hold for any sub-polytope of $\conv(B_{r,n})$. In 
this case, instead of the complete bipartite graph $K_{r,n-r}$, we would have a
subgraph $G$ of $K_{r,n-r}$. The maximal simplices of $\Delta_\succ$ would
correspond to maximal thrackles of $G$. It would be interesting to study the
number of such thrackles for other classes of matroids such as graphs,
transversals, etc.

It should be noted that the volume of any tangent cone for any matroid of rank
$r$ on $n$ elements is bounded by the volume of the tangent cone of the uniform
matroid $U^{r,n}$. That is, bounded by $n-2 \choose r-1$. Also immediate from
\autoref{thm:main} is that when the rank is fixed, the volume of the 
tangent cone is bounded polynomially in $n$. This provides an alternate proof 
of Lemma 10 in \cite{De-Loera:2009uq}.

Unfortunately, knowledge of the exact volume of the convex hull of a vertex of
a matroid polytope and its adjacent vertices does not immediately give a bound
on the volume of the matroid polytope itself. There are points in the matroid
polytope that are not in the convex hull of a vertex of a matroid polytope and
its adjacent vertices.


\bibliographystyle{plain}
\bibliography{references}

\begin{thebibliography}{10}

\bibitem{bar}
Alexander~I. Barvinok.
\newblock Polynomial time algorithm for counting integral points in polyhedra
  when the dimension is fixed.
\newblock {\em Mathematics of Operations Research}, 19:769--779, 1994.

\bibitem{De-Loera:2009uq}
Jesus~A. {De Loera}, David Haws, and Matthias K\"oppe.
\newblock Ehrhart polynomials of matroid polytopes and polymatroids.
\newblock {\em Journal of Discrete and Computational Geometry}, 2009.

\bibitem{Gelfand1987Combinatorial-g}
I.~M. Gelfand, M.~Goresky, R.~D. MacPherson, and V.~V. Serganova.
\newblock Combinatorial geometries, convex polyhedra, and {S}chubert cells.
\newblock {\em Adv. Math.}, 63:301--316, 1987.

\bibitem{Haws:2009uq}
David Haws.
\newblock Matroid polytopes: Algorithms, theory, and applications.
\newblock {\em arXiv:0905.4405}, 2009.

\bibitem{Loera:1995fk}
J.~De Loera, B.~Sturmfels, and R.~Thomas.
\newblock Gr\"obner bases and triangulations of the second hypersimplex.
\newblock {\em Combinatorica}, 15:409--424, 1995.

\bibitem{Oxley1992Matroid-Theory}
J.~Oxley.
\newblock {\em Matroid Theory}.
\newblock Oxford University Press, New York, NY, USA, 1992.

\bibitem{Stanley1996Combinatorics-a}
Richard~P. Stanley.
\newblock {\em Combinatorics and Commutative Algebra: Second Edition}.
\newblock Birkh\"auser, Boston, 2nd edition, 1996.

\bibitem{Stanley1997Enumerative-Com}
Richard~P. Stanley.
\newblock {\em Enumerative Combinatorics}, volume~1.
\newblock Cambridge University Press, 1997.

\bibitem{Sturmfels1996Grobner-Bases-a}
Bernd Sturmfels.
\newblock {\em Gr\"obner Bases and Convex Polytopes}, volume~8 of {\em
  University Lecture Series}.
\newblock American Mathematical Society, 1996.

\bibitem{Topkis1984Adjacency-on-Po}
D.~M. Topkis.
\newblock Adjacency on polymatroids.
\newblock {\em Mathematical Programming}, 30(2):229--237, October 1984.

\bibitem{Welsh1976Matroid-Theory}
D.~Welsh.
\newblock {\em Matroid Theory}.
\newblock Academic Press, Inc., 1976.

\end{thebibliography}

\section{appendix}

\begin{proof}[Proof of \autoref{thm:grobner}]
Note the reduction relation defined by the proposed Gr\"obner basis amounts to
replacing non-crossing edges by crossing edges. For each binomial $x_{ij}x_{kl}
- x_{il}x_{jk}$ in $\Cg$, the initial term with respect to $\succ$ corresponds
to the disjoint edges. This follows from the convex embedding of $K_{r,n-r}$
and the definition of the weights. The integral vectors in the kernel of
$\B_{r,n}$ are in bijection with even length closed walks on th complete
bipartite graph $K_{r,n-r}$, and hence so are the binomials of $I_{\B_{r,n}}$.
More precisely, with an even walk $\Gamma =
(i_1,i_2,\ldots,i_{2k-1},i_{2k},i_1)$ we associate the binomial

\begin{equation*}
b_\Gamma = \prod_{l=1}^k x_{i_{2l-1},i_{2l}} - \prod_{l=1}^k x_{i_{2l},i_{2l+1}}
\end{equation*}
where $i_{2k+1} = i_1$. Clearly the walk $\Gamma$ can be recovered from its
binomial $b_\Gamma$. By Corollary 4.4 of \cite{Sturmfels1996Grobner-Bases-a},
the infinite set of binomials associated with all even closed walks in
$K_{r,n-r}$ contains every reduced Gr\"obner basis of $I_{\B_{r,n}}$.
Therefore, in order to prove that $\Cg$ is a Gr\"obner basis, it is enough to
prove that the initial monomial of any binomial $b_\Gamma$ is divisible by some
monomial $x_{ij}x_{kl}$ where $(i,j),(k,l)$ is the pair of disjoint edges.

Suppose on the contrary there exists a binomial $b_\Gamma \in I_{\B_{r,n}}$
that contradicts our assertion. This implies that each pair of edges appearing
in the initial monomial of $b_\Gamma$ intersects. We may assume that $b_\Gamma$
is a minimial counter-example in the sense that $r,n$ are minimal and
$b_\Gamma$ has minimal weight. Hence the weight of the binomial $b_\Gamma$ is
the sum of the weights of its two terms. The walk $\Gamma$ is spanning in
$K_{r,n-r}$ by minimality of $n$. Every edge of $\Gamma$ gets a label "odd" or
"even" according to its position on the walk. In the case of the complete
bipartite graph, the odd edges go from $1 \leq i \leq r$ to $r+1 \leq j \leq
n$, and even edges go in the opposite direction.  If an edge is visited more
than once, it can not receive both "odd" and "even", since otherwise the
related variable can be factored out of $b_\Gamma$. This would contradict the
minimality of the weight. Alternatively, due to the fact above about odd and
even edges on the bipartite graph, this will not occur. Moreover, if $b_\Gamma
= \ve x^{\ve u} - \ve x^{\ve v}$ and $in_\succ(b_\Gamma) = \ve x^{\ve u}$, we
can assume that each pair of edges in $\ve x^{\ve v}$ intersects. Otherwise if
$(i,j),(k,l)$ is a non-intersecting pair of edges then we can reduce $\ve
x^{\ve v}$ modulo $\Cg$ to obtain a counterexample of smaller weight.

Suppose we draw $K_{r,n-r}$ as described in \autoref{rem:draw}. The 
\emph{circular distance} between any two vertices $1 \leq i \leq r$
and $r+1 \leq j \leq n$ is the shortest distance between $i$ and $j$
on the n-gon.

Let $(s,t)$ be an edge of the walk $\Gamma$ such that the circular distance
between $s$ and $t$ is smallest possible. The edge $(s,t)$ separates the
vertices of $K_{r,n-r}$, except $s$ and $t$, into two disjoint sets $P$ and $Q$
where $|P| \geq |Q|$. Let us assume the walk starts at $(s,t) = (i_1,i_2)$.
The walk is then a sequence of vertices and edges $\Gamma =
(i_1,(i_1,i_2),i_2,(i_2,i_3),\ldots,i_{2k},(i_{2k},i_{2k+1}))$.  Each pair of
odd (resp. even) edges intersects. The odd edges are of type
$(i_{2r-1},i_{2r})$ and the even edges are of type $(i_{2r},i_{2r+1})$.  Since
the circular distance of $(i_1,i_2)$ is minimal, the vertex $i_3$ can not be in
$Q$. Otherwise the edge $(i_2,i_3)$ would have smaller cicular distance. We
claim that if $P$ contains an odd vertex $i_{2r-1}$, then it also contains the
subsequent odd vertices $i_{2r+1},i_{2r+3},\ldots,i_{2k-1}$. The edge
$(i_1,i_2)$ is the common boundary of the two regions $P$ and $Q$. Any odd edge
intersects it (at least by having an end $\{i_1,i_2\}$) and thus $i_{2r}$ is in
$Q \cup \{i_1,i_2\}$. Since any even edge must intersect $(i_2,i_3)$, the
vertex $i_{2r+1}$ lies in $P \cup \{i_2\}$. To complete the proof of the claim
we show that $i_{2r+1} \neq i_2$. The equality $i_{2r+1} = i_2$ would imply
either $i_{2r}=i_1$ or $i_{2r} \in Q$. If $i_{2r} = i_1$ then $(i_1,i_2)$ is
both odd and even. On the other hand if $i_{2r} \in Q$ then $(i_{2r},i_2)$ has
smaller circular distance than $(i_1,i_2)$. Thus $i_{2r+1}$ belongs to $P$. The
claim is proved by repeating this argument.

Since $i_3$ was shown to be in $P$, it follows that all odd vertices except
$i_1$ lie in $P$ and the even vertices lie in $Q\cup \{i_1,i_2\}$. The final
vertex $i_{2k}$ is thus in $Q$. The even edge $(i_{2k},i_1)$ must be a closed
line segment contained in the region $Q$. Therefore $(i_1,i_3)$ and
$(i_{2k},i_1)$ are two even edges that do not intersect, which is a
contradiction. This proves $\Cg$ is a Gr\"obner basis of $I_{\B_{r,n}}$.

By construction, no monomial in an element of $\Cg$ is divisible by the initial
term of an element in $\Cg$. Hence $\Cg$ is the reduced Gr\"obner basis of
$I_{\B_{r,n}}$ with respect to $\succ$.
\end{proof}

\end{document}